\documentclass[11pt,letterpaper,reqno]{amsart}

\usepackage{amsmath}
\usepackage{amsfonts}
\usepackage{amssymb}
\usepackage{amsthm}
\usepackage{amscd}

\usepackage[colorlinks=true]{hyperref}
\hypersetup{urlcolor=green, citecolor=blue}

\usepackage{latexsym}
\usepackage{mathrsfs}
\usepackage{psfrag}
\usepackage{eufrak}

\usepackage[dvips]{epsfig}
\usepackage{epsfig}
\usepackage[all]{xy}         
\usepackage{wrapfig}

\usepackage[margin=2.5cm]{geometry}
\usepackage{marginnote}

\theoremstyle{plain}                              
\newtheorem{thm}{Theorem}[section]
\newtheorem{prop}{Proposition}
\newtheorem*{defn}{Definition}
\newtheorem{lem}[thm]{Lemma}
\newtheorem{cor}[thm]{Corollary}

\newtheorem*{thmA}{Theorem A}
\newtheorem*{thmB}{Theorem B}    

\newtheorem*{quest}{Question}

\theoremstyle{definition}                         
\newtheorem{example}{Example}
\newtheorem*{remark}{Remark} 

\theoremstyle{remark}                             

\numberwithin{equation}{section}

\newcommand{\R}{\mathbb{R}}                     

\newcommand{\ms}[1]{\mathscr{#1}}               

\newcommand{\veps}{\varepsilon}        

\newcommand{\del}{\partial}

\newcommand{\sn}{\text{\rm sn}}
\newcommand{\cn}{\text{\rm cn}}
\newcommand{\dn}{\text{\rm dn}}

\providecommand{\norm}[1]{\left\lVert#1\right\rVert}       

\providecommand{\abs}[1]{\left\lvert#1\right\rvert}        

\begin{document}

\title[Subriemannian geometry of contact Anosov flows]{On the
  subriemannian geometry of contact Anosov flows}


\author[S. N. Simi\'c]{Slobodan N. Simi\'c}

\address{Department of Mathematics and Statistics, San Jos\'e State University, San
  Jos\'e, CA 95192-0103}



\email{simic@math.sjsu.edu}



\date{\today} 





\begin{abstract}
  We investigate certain natural connections between subriemannian
  geometry and hyperbolic dynamical systems. In particular, we study
  dynamically defined horizontal distributions which split into two
  integrable ones and ask: how is the energy of a subriemannian
  geodesic shared between its projections onto the integrable
  summands? We show that if the horizontal distribution is the sum of
  the strong stable and strong unstable distributions of a special
  type of a contact Anosov flow in three dimensions, then for any
  short enough subriemannian geodesic connecting points on the same
  orbit of the Anosov flow, the energy of the geodesic is shared
  \emph{equally} between its projections onto the stable and unstable
  bundles. The proof relies on a connection between the geodesic
  equations and the harmonic oscillator equation, and its explicit
  solution by the Jacobi elliptic functions. Using a different idea,
  we prove an analogous result in higher dimensions for the geodesic
  flow of a closed Riemannian manifold of constant negative curvature.
\end{abstract}

\maketitle




\section{Introduction}
\label{sec:intro}

The goal of this paper is to investigate certain natural but insufficiently explored connections between hyperbolic dynamical systems and subriemannian (or Carnot-Carath\'eodory) geometry. A subriemannian geometry on a smooth connected manifold $M$ is a geometry defined by a nowhere integrable distribution $E$, called a \textsf{horizontal distribution}, equipped with a Riemannian metric $g$. Both $E$ and $g$ are required to be at least continuous but in most scenarios they are usually $C^\infty$. Since we can extend any partially defined Riemannian metric to the entire tangent bundle and the extension does not affect the properties of the subriemannian geometry, we will always assume that $g$ is defined on the entire tangent bundle. 

If $\gamma : [a,b] \to M$ is a horizontal (i.e., tangent to $E$) path, its length is defined in the usual way by
\begin{displaymath}
  \abs{\gamma} = \int_a^b \norm{\dot{\gamma}(t)} \: dt,
\end{displaymath}
where $\norm{v} = \sqrt{g(v,v)}$, for any vector $v \in E$.

A horizontal distribution $E$ on $M$ is called \textsf{nowhere integrable} if for every $p \in M$ and every $\veps > 0$ there exists a neighborhood $U$ of $p$ in $M$ such that every point in $U$ can be connected to $p$ by a horizontal path of length $< \veps$. In particular, every two points of $M$ can be connected by a horizontal path. This definition avoids certain undesirable pathological behavior which can arise if $E$ is not smooth; see \cite{sns+pams+10}.

The \textsf{subriemannian distance} between $x, y \in M$ is given by
\begin{displaymath}
  d_H(x,y) = \inf \{ \abs{\gamma} : \gamma \ \text{is a horizontal path from} \ x \ \text{to} \ y \}.
\end{displaymath}
A \textsf{subriemannian geodesic} from $x$ to $y$ is any horizontal path $\gamma$ which minimizes length among all horizontal paths connecting $x$ and $y$. Thus $\abs{\gamma} = d_H(x,y)$.

Recall that a $C^\infty$ horizontal distribution $E$ is called \textsf{bracket generating} if any local smooth frame $\{X_1, \ldots, X_k \}$ for $E$ together with all its iterated Lie brackets span the entire tangent bundle of $M$. (In the PDE literature, the bracket generating condition is called the H\"ormander condition.) By the Chow-Rashevskii theorem \cite{mont02} any bracket generating distribution is nowhere integrable.

\begin{figure}[htbp]
\centerline{
	\psfrag{x}[][]{$x$}
	\psfrag{y}[][]{$y$}
	\psfrag{c}[][]{$c$}
	\psfrag{g}[][]{$\gamma$}
        \psfrag{H}[][]{$E$}
        \psfrag{H1}[][]{$E_1$}
        \psfrag{H2}[][]{$E_2$}
        \psfrag{V}[][]{$V$}
\includegraphics[width=0.5\hsize]{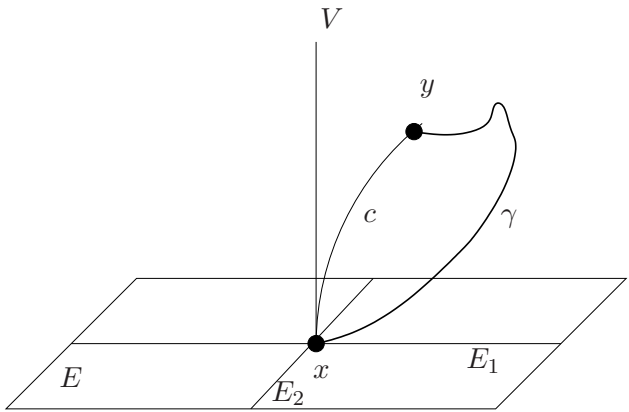}}
\caption{A subriemannian geodesic $\gamma$ connecting $x = c(0)$ and $y = c(1)$, where $c$ is a path tangent to the vertical bundle $V$.}
\label{fig:gen}
\end{figure}

It is sometimes the case that a horizontal distribution $E$ splits into two integrable orthogonal distributions, $E = E_1 \oplus E_2$, and $E$ is in turn orthogonal to a globally defined vertical distribution $V$, with $TM = E \oplus V$. If $E$ is bracket-generating, then any motion in the vertical direction is due to the fact that iterated Lie brackets of vector fields in $E_1$ and those in $E_2$ generate the entire tangent bundle. Given a ``vertical'' curve $c$ tangent to $V$ with endpoints $x$ and $y$ and a unit speed subriemannian geodesic $\gamma$ connecting $x$ and $y$, it is natural to ask the following question (see Figure~\ref{fig:gen}).

\begin{quest}
  If $TM = E \oplus V$ and $E = E_1 \oplus E_2$, how is the energy of a subriemannian geodesic $\gamma$ connecting endpoints of a curve tangent to $V$ shared between its projections onto $E_1$ and $E_2$?
\end{quest}
Stated more precisely, if $\gamma : [0,\ell] \to M$ is a unit speed horizontal path and $E = E_1 \oplus E_2$, then $\dot{\gamma}(t) = w_1(t) + w_2(t)$, with $w_i(t) \in E_i$. We define
\begin{displaymath}
  \mathcal{E}_i(\gamma) = \frac{1}{\ell} \int_0^\ell \norm{w_i(t)}^2 \: dt,
\end{displaymath}
for $i = 1, 2$, and think of $\mathcal{E}_i(\gamma)$ as the \emph{energy} of the projection of $\gamma$ to $E_i$. Clearly, $0 \leq \mathcal{E}_i(\gamma) \leq 1$ and $\mathcal{E}_1(\gamma) + \mathcal{E}_2(\gamma) = 1$. If $\mathcal{E}_1(\gamma) = \mathcal{E}_2(\gamma)$ we call $\gamma$ a $(E_1,E_2)$-\textsf{balanced horizontal path}. The above question therefore asks if every subriemannian geodesic connecting endpoints of a vertical path is $(E_1,E_2)$-balanced.

\begin{example}[The Heisenberg group]    \label{ex:heis}
  The Heisenberg group is a subriemannian geometry on $M = \R^3$ defined by the horizontal distribution $E$ which is the kernel of the 1-form $\alpha = dz - \frac{1}{2}(x \: dy - y \: dx)$. The Riemannian metric on $E$ is defined by $ds^2 = dx^2 + dy^2$. The vector fields
  \begin{displaymath}
    X_1 = \frac{\del}{\del x} - \frac{y}{2} \frac{\del}{\del z} \quad \text{and} \quad
    X_2 = \frac{\del}{\del y} - \frac{x}{2} \frac{\del}{\del z}
  \end{displaymath}
  form a global orthonormal frame for $E$. It is not hard to check that $[X_1,X_2] = \frac{\del}{\del z} =: X_0$, so $E$ is bracket-generating. Since $[X_0, X_1] = [X_0,X_2] = 0$, the Heisenberg group is nilpotent. We will show in \S\ref{sec:sr-geom} that every Heisenberg subriemannian geodesic whose endpoints differ only in the $z$-coordinate is $(E_1,E_2)$-balanced. This follows easily from the fact that Heisenberg geodesics are lifts of circles in the $xy$-plane.
\end{example}

Subriemannian geometries whose horizontal distributions have a natural splitting into two integrable distributions occur frequently in hyperbolic and partially hyperbolic dynamical systems. For instance, if $f$ is a partially hyperbolic diffeomorphism of a compact manifold $M$, then $f$ preserves two invariant bundles called the stable $E^s$ and the unstable $E^u$ bundles, both uniquely integrable. Transverse to them is the center bundle $E^c$, which is not always integrable. Although $E^s$ and $E^u$ are usually not smooth, their sum $E^{su} = E^s \oplus E^u$ frequently has the so called \textsf{accessibility property}. This means that any two points in $M$ can be joined by a continuous piecewise smooth path whose smooth legs are alternately tangent to $E^s$ and $E^u$. Thus $E^{su}$ naturally defines a subriemannian geometry on $M$ and we can take $H = E^{su}$ and $V = E^c$. Accessibility plays an important role in partially hyperbolic dynamics where it is an essential ingredient in the study of stably ergodic systems and the Pugh-Shub conjecture \cite{ps+04}. The main difficulty with the subriemannian geometry defined by $E^{su}$ is that it lacks smoothness, so it is not amenable to analysis using standard techniques.

In this paper we consider the case of \emph{contact Anosov flows}, where the natural horizontal distribution is always at least $C^1$. This is a scenario which is in a sense diametrically opposite to that of the Heisenberg group.

Recall that a non-singular smooth flow $\Phi = \{f_t \}$ on a closed (compact and without boundary) Riemannian manifold $M$ is called an \textsf{Anosov flow} if there exists an invariant splitting $TM = E^{ss} \oplus E^c \oplus E^{uu}$ such that $E^c$ is spanned by the infinitesimal generator $X$ of the flow, $E^{ss}$ is uniformly exponentially contracted and $E^{uu}$ is uniformly exponentially expanded by the flow in positive time. We call $E^{ss}$ and $E^{uu}$ the \textsf{strong stable} and \textsf{strong unstable bundles}; $E^c$ is the \textsf{center bundle}.

A \textsf{contact structure} on a manifold $M$ of dimension $2n+1$ is a $C^1$ hyperplane field $E$ which is as far from being integrable as possible~\cite{duff+solomon+sympl+99}. This means that there exists a $C^1$ 1-form $\alpha$ such that $\text{Ker}(\alpha) = E$ and $\alpha \wedge (d\alpha)^n$ is a volume form for $M$; $\alpha$ is called a \textsf{contact form} for $E$. Contact structures are always bracket-generating.

A vector field $X$ is called the \textsf{Reeb vector field} of $\alpha$ if $\alpha(X) = 1$ and $X$ is in the kernel of $d\alpha$, i.e., $i_X d\alpha = 0$. An Anosov flow is called \textsf{contact} if $E = E^{su}$ is a contact structure (and in particular $C^1$) and the infinitesimal generator $X$ of the flow is the Reeb vector field for the contact form $\alpha$ for $E^{su}$ with $\alpha(X) = 1$. Our goal is to understand the subriemannian geometry defined by the distribution $E^{su}$ associated with a contact Anosov flow. We will call subriemannian geodesics of this geometry $su$-\textsf{subriemannian geodesics}. An $su$-subriemannian geodesic will be called $su$-\textsf{balanced} if it is balanced with respect to the splitting $E^{ss} \oplus E^{uu}$.

Contact Anosov flows have good dynamical properties; in particular, they exhibit exponential decay of correlations (cf., \cite{liverani+contact+04}). Until recently however, the only known contact Anosov flows were the geodesic flows of Riemannian or Finsler manifolds; in \cite{foulon+hass+13} Foulon and Hasselblatt used surgery near a transverse Legendrian knot to construct many new contact Anosov flows on 3-manifolds which are not topologically orbit equivalent to any algebraic flow.

Assume now that $\Phi$ is a contact Anosov flow on a 3-manifold $M$. Denote its infinitesimal generator by $X$ and let $Y$ and $Z$ be unit (with respect to some Riemannian metric $g$ whose volume form equals the contact volume form) vector fields in $E^{ss}$ and $E^{uu}$ respectively. Then $T_xf_t(Y) = \mu(x,t) Y$ and $T_xf_t(Z) = \lambda(x,t) Z$, for some 1-cocycles $\mu, \lambda : M \times \R \to \R$, with $\lambda \mu = 1$, where $Tf_t$ denotes the tangent map (i.e., derivative) of the time-$t$ map $f_t$ of the flow. Thus
\begin{equation}    \label{eq:bracket}
  [X,Y] = aY \qquad \text{and} \qquad [X,Z] = -aZ,
\end{equation}
where $a(x) = -\dot{\mu}(x,0) = \dot{\lambda}(x,0)$, for all $x \in M$. Since $E^{su}$ is contact, it follows that $[Y,Z]$ is transverse to $E^{su}$.

\begin{defn}
  An Anosov flow $\Phi$ on a 3-dimensional closed manifold will be called a \textsf{special contact Anosov flow} if there exists a $C^1$ Riemannian metric $g$ and a $C^1$ global orthonormal frame $(X,Y,Z)$ relative to $g$ such that:

  \begin{enumerate}
  \item[(a)] $E^c = \R X$, $E^{ss} = \R Y$ and $E^{uu} = \R Z$.

  \item[(b)] $[X,Y] = Y$, $[Y,Z] = X$ and $[Z,X] = Z$.

  \end{enumerate}
\end{defn}

Note that $X, Y, Z$, and $g$ are required to be only $C^1$. The following lemma shows that must in fact be $C^\infty$.

\begin{lem}   \label{lem:smooth}
  If $\Phi$ is a special contact Anosov flow, then (with the notation as above), $X, Y$ and
 $Z$ are all $C^\infty$.
\end{lem}

\begin{proof}
  Let $(\alpha,\beta,\gamma)$ be the coframe dual to $(X,Y,Z)$. Since $X, Y$ and $Z$ are $C^1$, so are $\alpha, \beta$ and $\gamma$. We have:
  \begin{align*}
    1 & = \alpha(X) \\
    & = \alpha([Y,Z]) \\
    & = Y \alpha(Z) - Z \alpha(Y) - d\alpha(Y,Z) \\
    & = - d\alpha(Y,Z).
  \end{align*}
  We can show in a similar way that $d\alpha(X,Y) = d\alpha(X,Z) = 0$. Therefore, $d\alpha$ is $C^1$ relative to a $C^1$ frame, hence $\alpha$ is, in fact, $C^2$. It follows analogously that $\beta$ and $\gamma$ are also $C^2$. Hence $X, Y$ and $Z$ are all $C^2$ as well. By bootstrap, it follows that $X, Y$ and $Z$ are in fact $C^\infty$.
\end{proof}

\begin{remark}
  Observe that if $\Phi$ is a special contact Anosov flow, then $X, Y$ and $Z$ span a copy of the Lie algebra $\mathfrak{sl}(2,\R)$, so the universal cover of $M$ is the Lie group $\text{SL}(2,\R)$ and the lift of the Anosov flow to the universal cover is an algebraic one. In particular, $M$ is a quotient of $\text{SL}(2,\R)$ by a discrete cocompact subgroup. In other words, special contact Anosov flows are precisely algebraic Anosov flows.

  Ghys~\cite{ghys+87} showed that in three dimensions every contact Anosov flow $\Phi$ with $C^\infty$ strong bundles is $C^\infty$ equivalent to an algebraic flow on a quotient $N = \Gamma\setminus\widetilde{\text{SL}}(2,\R)$, in the sense that there exists a $C^\infty$ diffeomorphism $h : N \to M$ that sends the orbits of the ``diagonal'' flow on $N$ to the orbits of $\Phi$. Therefore, every contact Anosov flow with $C^\infty$ strong bundles is $C^\infty$ orbit equivalent to a special Anosov flow.
\end{remark}

Our main results are the following.

\begin{thmA}
  Let $\Phi = \{ f_t \}$ be a special contact Anosov flow on a closed Riemannian 3-manifold $M$. Then there exists a $\delta > 0$ such that for all $x \in M$ and $\abs{t} < \delta$, every $su$-subriemannian geodesic connecting $x$ and $f_t(x)$ is $su$-balanced.
\end{thmA}

In higher dimensions we prove a result analogous to Theorem A if the contact Anosov flow is the geodesic flow on the unit tangent bundle of a manifold with \emph{constant} negative sectional curvature.

\begin{thmB}
  Let $\Phi = \{ f_t \}$ be the geodesic flow of a closed Riemannian manifold $N$ with
  constant negative sectional curvature on its unit tangent bundle $M$. Then there exists a $\delta
  > 0$ such that for all $x \in M$ and $\abs{t} < \delta$, every $su$-subriemannian geodesic connecting $x$ and $f_t(x)$ is $su$-balanced.
\end{thmB}

\paragraph{\textbf{Outline of the paper.}} In Section \S\ref{sec:prelim} we review some basic results on Anosov and geodesic flows, subriemannian geodesics, and the solution of the harmonic oscillator equation via Jacobi elliptic functions. Theorem A is proved in Section \S\ref{sec:proofA} and Theorem B in Section \S\ref{sec:proofB}. We conclude the paper with a list of open questions in Section \S\ref{sec:questions}.

\subsection*{Acknowledgments}

We are grateful to Alan Weinstein, who generously offered the main idea of proof of Theorem B. We would also like to thank him for many inspiring conversations and moral support over the years.

\section{Preliminaries}
\label{sec:prelim}

\subsection{Anosov flows}
\label{sec:anosov}

A non-singular smooth flow $\Phi = \{f_t \}$ on a closed Riemannian
manifold $M$ is called an \textsf{Anosov flow} if there exists an
invariant splitting $TM = E^{ss} \oplus E^c \oplus E^{uu}$ such that
$E^c$ is spanned by the infinitesimal generator of the flow and there
exist uniform constants $c > 0$, $0 < \mu_- \leq \mu_+ < 1$ and
$\lambda_+ \geq \lambda_- > 1$ such that for all $v \in E^{ss}$, $w
\in E^{uu}$, and $t \geq 0$, we have
\begin{equation}      \label{eq:ss}
  \frac{1}{c} \mu_-^t \norm{v} \leq \norm{Tf_t (v)} \leq c \mu_+^t \norm{v},
\end{equation}
and
\begin{equation}   \label{eq:uu}
  \frac{1}{c} \lambda_-^t \norm{w} \leq \norm{Tf_t (w)} \leq c \lambda_+^t \norm{w}.
\end{equation}
The strong stable $E^{ss}$ and strong unstable bundles $E^{uu}$ are in
general only H\"older continuous \cite{hps77}, but they are
nevertheless always uniquely integrable giving rise to the strong
stable and strong unstable foliations denoted by $W^{ss}$ and
$W^{uu}$, respectively. The codimension one distribution $E^{su} =
E^{ss} \oplus E^{uu}$ is generally not integrable; if it is, then by
Plante~\cite{plante72}, the flow admits a global cross section and is
therefore topologically conjugate to a suspension of an Anosov
diffeomorphism.

The bundles $E^{cs} = E^c \oplus E^{ss}$ and $E^{cu} = E^c \oplus
E^{uu}$ are called the \textsf{center stable} and \textsf{center
  unstable bundles}. They are generically only H\"older continuous
\cite{hps77}, but are always uniquely integrable
\cite{anosov+67}. However, if $\dim M = 3$ and the flow is $C^3$ and
preserves the Riemannian volume, then it follows from the work of
Hurder and Katok \cite{hurder+katok+90} that $E^{cs}$ and $E^{cu}$ are
both of class $C^1$ and the transverse derivatives of both bundles are
$C^\theta$-H\"older, for \emph{all} $0 < \theta < 1$.

Without loss we will always assume that all the invariant bundles of
an Anosov flow are orientable. (If not, we can pass to a double cover.)

We will need the following easy lemma.

\begin{lem}       \label{lem:strong}
  Let $\Phi$ be a contact Anosov flow on a 3-manifold $M$. Then
  $E^{uu}$ and $E^{ss}$ are both $C^1$.
\end{lem}

\begin{proof}
  Since $E^{su}$ is $C^1$ by assumption and $E^{cs}$ and $E^{cu}$ are
  $C^1$ by \cite{hurder+katok+90}, it follows that $E^{ss} = E^{cs}
  \cap E^{su}$ and $E^{uu} = E^{cu} \cap E^{su}$ are also $C^1$.
\end{proof}

\subsection*{Geodesic flows}

In this section we briefly review some basic facts about geodesic
flows. If $N$ is a Riemannian manifold, then its geodesic flow $\Phi =
\{ f_t \}$ restricted to the unit tangent bundle $M = T^1 N$ of $N$
admits a canonical contact form (cf.,
\cite{paternain+geodesic+99}). If the sectional curvature $K$ of $N$
is negative, then $\Phi$ is known to be of Anosov type
\cite{anosov+67,eberlein+73}, in which case $E^{su}$ is a contact
structure and $TM = E^c \oplus E^{su}$ is an orthogonal splitting with
respect to the Sasaki metric \cite{paternain+geodesic+99}. If the
sectional curvature $K$ is constant (and negative), then $E^{ss}$ and
$E^{uu}$ are $C^\infty$, but if $K$ is variable, then $E^{ss}$ and
$E^{uu}$ are only of class $C^{1+\theta}$, for some $0 < \theta <
1$~\cite{hp+75}.

Assume now that $K$ is constant and negative. Without loss we can
assume that $K = -1$. Then (cf., \cite{anosov+67}) there exists a
constant $c > 0$ such that
\begin{displaymath}
  \norm{Tf_t(v)} = e^{-ct} \norm{v} \quad \text{and} \quad
  \norm{Tf_t(w)} = e^{ct} \norm{w},
\end{displaymath}
for all $t \in \R$, $v \in E^{ss}$ and $w \in E^{uu}$. In other words,
the flow contracts all stable directions and expands all unstable
directions at the same rates at all points of $M$. For simplicity, we
will assume that $c = 1$; this can always be achieved by a constant
time change.

Let $F$ be an isometry of $N$ and denote by $F_\ast$ the restriction
of $TF$ to the unit tangent bundle $M = T^1 N$ of $N$. Then for any
unit-speed geodesic $t \mapsto c(t)$ in $N$, we have
$F_\ast(\dot{c}(t)) = F_\ast(f_t \dot{c}(0)) = f_t
(F_\ast(\dot{c}(0))$, which implies that $F_\ast$ preserves the
geodesic vector field $X$. We claim that $F_\ast$ also preserves the
strong stable $W^{ss}$ and strong unstable $W^{uu}$ foliations of the
geodesic flow. Indeed, since $F_\ast \circ f_t = f_t \circ F_\ast$,
for any $v_1, v_2$ in the same $W^{ss}$-leaf, we have:
\begin{align*}
  d(f_t(F_\ast(v_1)), f_t(F_\ast(v_2))) & = d (F_\ast(f_t(v_1)),
  F_\ast(f_t(v_2))) \\
  & \leq K d(f_t(v_1),f_t(v_2)) \\
  & \to 0,
\end{align*}
as $t \to \infty$, where $K$ is the Lipschitz constant of $F_\ast$
(which is finite, since $M$ is compact and $F_\ast$ is smooth) and $d$
denotes the distance function on $M$ induced by the Sasaki
metric. Therefore, $F_\ast(v_1)$ and $F_\ast(v_2)$ lie in the same
$W^{ss}$-leaf. Thus $F_\ast$ preserves $W^{ss}$ and $TF_\ast$
preserves $E^{ss}$. It can similarly be shown that $E^{uu}$ is also
invariant with respect to $TF_\ast$. Thus $TF_\ast$ preserves the
splitting $E^c \oplus E^{su}$. An analogous statement is true for any
lift $\tilde{F}_\ast$ of $F_\ast$ to the universal Riemannian covering
space $\tilde{M}$ of $M$.

Recall that if $K = -1$, then $N = \mathbb{H}^n/\Gamma$, where
$\Gamma$ is a group of isometries of $\mathbb{H}^n$ acting freely and
properly discontinuously on it~\cite{boothby+2003}. It is clear that
$\Gamma$ also acts freely and properly discontinuously on the unit
tangent bundle $T^1 \mathbb{H}^n$ of $\mathbb{H}^n$ and that $(T^1
\mathbb{H}^n)/\Gamma$ is isometric to $T^1(\mathbb{H}^n/\Gamma) = T^1
N = M$. Thus $\tilde{M}$ is isometric to the universal Riemannian
covering space of $(T^1\mathbb{H}^n)/\Gamma$. Since $n \geq 3$, $T^1
\mathbb{H}^n$ is simply connected, so it \emph{is} the universal
covering space of $(T^1\mathbb{H}^n)/\Gamma$. Thus $\tilde{M}$ is
isometric to $T^1 \mathbb{H}^n$.

\begin{lem}    \label{lem:hyperbolic}
  For all $\tilde{u}, \tilde{v} \in \tilde{M}$ there exists an
  isometry $F$ of $\tilde{M}$ such that $F(\tilde{u}) = \tilde{v}$ and
  $F$ leaves the lift $\tilde{X}$ of the geodesic vector field $X$
  invariant.
\end{lem}

\begin{proof}
  By the above observation, $\tilde{M}$ is isometric to $T^1
  \mathbb{H}^n$, so we can identify $\tilde{u}, \tilde{v}$ with unit
  tangent vectors to $\mathbb{H}^n$ at some points $x, y \in
  \mathbb{H}^n$, respectively. There exists an isometry $f$ of
  $\mathbb{H}^n$ such that $f(x) = y$ and $Tf(\tilde{u}) = \tilde{v}$
  (see \cite{boothby+2003}). Since isometries map geodesics to
  geodesics, $Tf$ leaves $\tilde{X}$ invariant. Thus $F =
  Tf \! \restriction_{T^1 \mathbb{H}^n}$ has the desired properties.
\end{proof}

\subsection{Subriemannian geodesics}
\label{sec:sr-geom}

In this section we briefly review subriemannian geodesic equations. We
follow \cite{mont02}.

Let $E$ be a bracket-generating distribution on a smooth manifold
$M$. For each smooth vector field $X$ on $M$ we define the
\textsf{momentum function} $P_X : T^\ast M \to \R$ by
\begin{displaymath}
  P_X(p) = p(X),
\end{displaymath}
for any $p \in T^\ast M$, where $T^\ast M$ is the cotangent bundle of
$M$. Thus the momentum function of $X$ is just the evaluation of any
covector on $M$ at $X$.

The \textsf{subriemannian Hamiltonian} $H$ of $E$ is the map $H : T^\ast M \to
\R$ defined by 
\begin{displaymath}
  H(p) = \frac{1}{2} \langle p, p \rangle,
\end{displaymath}
where $\langle \cdot, \cdot \rangle$ denotes the cometric on $T^\ast
M$ induced by the Riemannian metric $g$ on $E$ (see \cite{mont02}). If
$(X_1,\ldots, X_k)$ is a local horizontal frame and $g_{ij} =
g(X_i,X_j)$, then the Hamiltonian can be expressed as
\begin{displaymath}
  H = \frac{1}{2} \sum_{i,j} g^{ij} P_{X_i} P_{X_j},
\end{displaymath}
where $g^{ij}$ are the entries of the inverse of the matrix
$[g_{ij}]$. In particular, if $(X_1,\ldots,X_k)$ is a local
\emph{orthonormal} frame for $E$, then
\begin{displaymath}
  H = \frac{1}{2} \sum_{i=1}^k P_{X_i}^2.
\end{displaymath}
The \textsf{normal geodesic equation} for $E$ is the equation
\begin{equation}       \label{eq:geod}
  \dot{f} = \{ f,H \},
\end{equation}
where $f : T^\ast M \to \R$ is a smooth function and $\{f,H \}$
denotes the Poisson bracket of $f$ and $H$. Projections of the
solutions to \eqref{eq:geod} to $M$ are called \textsf{normal
  geodesics}.

Recall that 
\begin{displaymath}   
  \{ f,H \} = \omega(X_f,X_H),
\end{displaymath}
where $\omega$ is the (canonical) symplectic form on $T^\ast M$ and
$X_f, X_H$ are the Hamiltonian vector fields defined by $f, H$,
respectively. It is well-known that the Poisson bracket defines a Lie
algebra structure on the ring of smooth functions on $T^\ast M$ and
that the map $f \mapsto \{ f, H \}$ satisfies the Leibniz rule thus
defining a vector field on $T^\ast M$ (which of course is exactly
$X_H$). Recall also that $\{ P_X, P_Y \} = - P_{[X,Y]}$, for any smooth
vector fields $X, Y$ on $M$.

The equation \eqref{eq:geod} is to be interpreted in the following way: if
$t \mapsto p(t)$ is an integral curve of the Hamiltonian vector field
$X_H$ and if $f : T^\ast M \to \R$ is any smooth function, then
\begin{displaymath}
  \frac{d}{dt} f(p(t)) = \{ f,H \}(p(t)).
\end{displaymath}
In canonical coordinates $(x_1,\ldots, x_n; p_1,\ldots,p_n)$, where
$(x_1,\ldots,x_n)$ are local coordinates on $M$ and $p_i =
P_{\del/\del x_i}$, the geodesics equations assume the familiar form:
\begin{displaymath}
  \dot{x}_i = \frac{\del H}{\del p_i}, \qquad \dot{p}_i = -\frac{\del
   H}{\del x_i}.
\end{displaymath}
\begin{thm}[\cite{mont02}]
  Let $t \mapsto \Gamma(t)$ be a solution to the normal geodesic
  equation \eqref{eq:geod} and let $\gamma$ be its projection to
  $M$. Then every sufficiently short arc of $\gamma$ is a
  subriemannian geodesic. If $E$ is a 2-step distribution, then every
  subriemannian geodesic is normal.
\end{thm}

Recall that $E$ is a 2-step distribution if for any local frame
$X_1,\ldots, X_k$ for $E$, the vector fields $X_1,\ldots,X_k$ together
with their first-order Lie brackets $[X_i,X_j]$ $(1 \leq i, j \leq k)$
generate the entire tangent bundle.

Assume now $\dim M = 3$ and $E$ is a contact structure. It is easy to
see that $E$ is a 2-step distribution. Let $(X_1,X_2)$ be a local
orthonormal frame for $E$ and $\alpha$ a contact form for $E$. Denote
the Reeb field of $\alpha$ by $X_0$. Clearly, $(X_0,X_1,X_2)$ is a
local frame for $TM$. The structure contants of the frame
$(X_0,X_1,X_2)$ are smooth functions $c_{ij}^k$ defined by
\begin{displaymath}
  [X_i,X_j] = \sum_{k=0}^2 c_{ij}^k X_k.
\end{displaymath}
It follows that
\begin{displaymath}
  \{ P_{X_i}, P_{X_j} \} = - \sum_{k=0}^2 c_{ij}^k P_{X_k}.
\end{displaymath}
The subriemannian Hamiltonian corresponding to the frame $(X_1,X_2)$
is 
\begin{displaymath}
  H = \frac{1}{2} (P_{X_1}^2 + P_{X_2}^2).
\end{displaymath}
Introduce fiberwise coordinates $(P_{X_0},P_{X_1},P_{X_2})$ on $T^\ast
M$. In these coordinates the normal geodesic equations are
\begin{align*}
  \dot{x} & = P_{X_1} X_1 + P_{X_2} X_2 \\
  \dot{P}_{X_i} & = \{ P_{X_i}, H \}, 
\end{align*}
for $i = 0, 1, 2$. 
\begin{example}[The Heisenberg group, continued]
  We will show that every Heisenberg subriemannian geodesic whose
  endpoints differ only in the $z$-component is balanced with respect
  to the splitting $E = E_1 \oplus E_2$. See
  Example~\ref{ex:heis}. Since $(X_1,X_2)$ is an orthonormal frame for
  $E$, the subriemannian Hamiltonian is
  \begin{displaymath}
    H = \frac{1}{2} (P_{X_1}^2 + P_{X_2}^2)
  \end{displaymath}
  and the subriemannian geodesic equation is $\dot{f} = \{ f, H \}$.

  Using $[X_1,X_2] = \del/\del z =: X_0$, $[X_0,X_1] = [X_0,X_2] = 0$,
  we obtain $\{ P_{X_1},P_{X_2} \} = -P_{X_0}$, $\{ P_{X_0},P_{X_1} \} =
  \{ P_{X_0},P_{X_2} \} = 0$. Therefore the subriemannian geodesic
  equations are
  \begin{align*}
    \dot{p} & = P_{X_1} X_1 + P_{X_2} X_2 \\
    \dot{P}_{X_0} & = 0 \\
    \dot{P}_{X_1} & = -P_{X_0} P_{X_2} \\
    \dot{P}_{X_2} & = P_{X_0} P_{X_1},
  \end{align*}
  where $p = (x,y,z)$. Since geodesics travel at constant speed, we
  can restrict the equations to the level set $P_{X_1}^2 + P_{X_2}^2 =
  1$ of $H$ and reparametrize $P_{X_1}$ and $P_{X_2}$ by
  \begin{displaymath}
    P_{X_1} = \cos \theta, \quad P_{X_2} = \sin \theta.
  \end{displaymath}
  It is not hard to check that the last three geodesic equations are
  equivalent to
  \begin{displaymath}
    \dot{\theta} = P_{X_0}, \qquad \dot{P}_{X_0} = 0.
  \end{displaymath}
  Thus $\ddot{\theta} = 0$, so $\theta(t) = v_0 t + \theta_0$, where $v_0
  = \dot{\theta}(0) = P_{X_0}(0)$ and $\theta_0 = \theta(0)$. It
  follows that every Heisenberg geodesic satisfies 
  \begin{displaymath}
    \dot{p} = \cos(v_0 t + \theta_0) X_1 + \sin(v_0 t + v_0) X_2,
  \end{displaymath}
  with real parameters $v_0$ and $\theta_0$ as above.  Note that
  $\dot{x} = \cos(v_0 t + \theta_0)$ and $\dot{y} = \sin(v_0 t + \theta_0)$.

  Now assume that a subriemannian geodesic $\gamma : [0,\ell] \to
  \R^3$ connects two points which differ only in the $z$-coordinate,
  i.e., they lie on an orbit of the flow of $X_0$. Projecting to the
  $xy$-plane we obtain
  \begin{displaymath}
    \int_0^\ell \cos(v_0 t + \theta_0) \: dt = x(\ell) - x(0) = 0, 
    \qquad \int_0^\ell \sin(v_0 t +
    \theta_0) \: dt = y(\ell) - y(0) = 0.
  \end{displaymath}
  Thus $v_0 \ell$ must be an integer multiple of $2\pi$. Let us show
  that $\gamma$ is balanced with respect to the splitting $E = \R X_1
  \oplus \R X_2$. We have
  \begin{displaymath}
    \mathcal{E}_1(\gamma) = \frac{1}{\ell} \int_0^\ell \cos^2(v_0 t + \theta_0) \:
    dt \quad \text{and} \quad 
    \mathcal{E}_2(\gamma) = \frac{1}{\ell} \int_0^\ell \sin^2(v_0 t + \theta_0) \:
    dt.
  \end{displaymath}
  Hence
  \begin{align*}
    \mathcal{E}_1(\gamma) - \mathcal{E}_2(\gamma) & = \frac{1}{\ell} \int_0^\ell \{
    \cos^2(v_0 t + \theta_0) - \sin^2(v_0 t + \theta_0) \} \: dt \\
    & = \frac{1}{\ell} \int_0^\ell \cos 2(v_0 t+\theta_0) \: dt \\
    & = 0,
  \end{align*}
  since $v_0 \ell = 2\pi n$, for some integer $n$. 
\end{example}

\subsection{Harmonic oscillator and Jacobi elliptic functions}
\label{sec:jacobi}

To make the paper as self-contained as possible, we review in some
detail the method of explicitly solving the harmonic oscillator (i.e.,
unforced undamped pendulum) equation
\begin{equation}   \label{eq:pendulum}
  \ddot{\theta} + \omega^2 \sin \theta = 0
\end{equation}
by the Jacobi elliptic functions $\sn$ and $\cn$, defined below. We
closely follow \cite{meyer+amm+01}, adding results we need along the way.

Let $0 < k < 1$. The \textsf{Jacobi elliptic functions} $\sn(t,k)$, $\cn(t,k)$
and $\dn(t,k)$ are defined as the unique solutions $x(t), y(t)$ and
$z(t)$ of the system of differential equations
\begin{align*}
  \dot{x} & = yz \\
  \dot{y} & = - zx \\
  \dot{z} & = -k^2 xy,
\end{align*}
satisfying the initial conditions
\begin{displaymath}
  x(0) = 0, \quad y(0) = 1, \quad z(0) = 1.
\end{displaymath}
The paramater $k$ is called the modulus. Some basic properties of $\sn, \cn$
and $\dn$ are listed in the following proposition whose proof can be
found in \cite{meyer+amm+01}.

\begin{prop}   \label{prop:elliptic}
  \begin{enumerate}

  \item[(a)] The Jacobi elliptic functions $\sn$, $\cn$ and $\dn$ are
    analytic and defined for all real $t$.

  \item[(b)] $\sn^2(t,k) + \cn^2(t,k) = 1$ and $k^2 \sn^2(t,k) +
    dn^2(t,k) = 1$, for all $t \in \R$ and $0 < k < 1$.

  \item[(c)] Let $K = K(k) > 0$ be the unique number such that
    $\cn(K,k) = 0$ and $\cn(t,k) > 0$, for all $0 < t < K$. That is,
    $K$ is the time it takes $\cn(t,k)$ to decrease to 0 from its
    initial value 1. Then $\sn(t,K)$ and $\dn(t,k)$ are even about $K$
    and $\cn(t,k)$ is odd about $K$.

  \item[(d)] $\sn(t,k)$ and $\cn(t,k)$ are $4K$-periodic in $t$ and
    $\dn(t,k)$ is $2K$-periodic in $t$.

  \item[(e)] The function $x(t) = \sn(t,k)$ is the unique solution to
    the initial value problem
    \begin{displaymath}
      (\dot{x})^2 = (1-x^2) (1-k^2 x^2), \qquad x(0) = 0, \quad
      \dot{x}(0) = 1.
    \end{displaymath}

  \end{enumerate}
\end{prop}

Observe that $\sn(t,k)$ and $\cn(t,k)$ have the same symmetries with
respect to $K$ as $\sin t$ and $\cos t$ have with respect to $\pi/2$.

\begin{cor}  \label{cor:elliptic}
  If 
    \begin{displaymath}
      \int_a^b \sn(t;k) \: dt = \int_a^b \cn (t;k) \: dt = 0,
    \end{displaymath}
    then $b-a$ is an integer multiple of $4K$. Moreover,
    \begin{displaymath}
      \int_a^b \sn(t;k) \: \cn(t;k) \: dt  = 0.
    \end{displaymath}
\end{cor}

\begin{proof}
  The proof follows from parts (c) and (d) of the previous
  Proposition. The calculations are analogous to those
  proving similar properties for the functions $\sin$ and $\cos$.
\end{proof}

Now consider the pendulum equation $\ddot{\theta} + \omega^2 \sin
\theta = 0$, with $\omega > 0$. It is not hard to check that the
``energy'' of the oscillator given by
\begin{displaymath}
  I = \frac{1}{4} \dot{\theta}^2 + \frac{\omega^2}{2} (1-\cos \theta)
\end{displaymath}
is constant along solutions. We can rewrite $I$ as
\begin{displaymath}
  I = \frac{1}{4} \dot{\theta}^2 + \omega^2 \sin^2 \frac{\theta}{2}.
\end{displaymath}
Let $y(t) = \sin \dfrac{\theta(t)}{2}$. Then $\dot{y} = \dfrac{1}{2} \dot{\theta}
\cos \dfrac{\theta}{2} = \pm \dot{\theta} \sqrt{1-y^2}$. Squaring both
sides and solving for $\dot{\theta}^2$ from the equation for $I$, we
obtain
\begin{equation}     \label{eq:y}
  (\dot{y})^2 = (1-y^2) (I-\omega^2 y^2).
\end{equation}
There are four possibilities.

\begin{description}

\item[Case 1] $I=0$. Then $\theta(t) \equiv 0$ (mod $2\pi$) and the
  pendulum is in the stable downward equilibrium.

\vspace{3mm}

\item[Case 2] $0 < I < \omega^2$. We look for a solution in the form
  $y(t) = A \: \sn(B(t-t_0);k)$, for some constants $A, B, t_0$ and $0
  < k < 1$, and obtain
  \begin{displaymath}
    y(t) = k \: \sn(\omega(t-t_0);k), \qquad \text{with} \quad k = \frac{\sqrt{I}}{\omega}.
  \end{displaymath}
Therefore, for any $t_0 \in \R$,
\begin{displaymath}
  \theta(t) = 2 \arcsin \left\{ k \:   
    \sn(\omega(t-t_0);k) \right\},
\end{displaymath}
is a solution to \eqref{eq:pendulum}. This case corresponds to the
pendulum swinging back and forth. Given a particular solution
$\theta(t)$, we can compute $t_0$ using $\theta(t_0) = 0$.

\vspace{3mm}

\item[Case 3] $I = \omega^2$. If $\dot{\theta}(0) = 0$, then
  $\theta(t) \equiv \pi$ (mod $2\pi$) and the pendulum is in the
  unstable upward equilibrium. If $\dot{\theta}(0) \neq 0$, then
  \begin{displaymath}
    \dot{y} = \pm \omega (1-y^2),
  \end{displaymath}
  and it is not hard to check that
  \begin{displaymath}
    y(t) = \pm \tanh(\omega(t-t_0))
  \end{displaymath}
  satisfies \eqref{eq:y} for any $t_0$, so the solution to the
  pendulum equation in this case is
  \begin{displaymath}
    \theta(t) = \pm 2 \arcsin \tanh(\omega(t-t_0)).
  \end{displaymath}
  This solution is a saddle connection connecting two downward equilibria.

\vspace{2mm}

\item[Case 4] $I > \omega^2$. We seek a solution in the form $y(t) =
  \sn(A(t-t_0);k)$, for some constants $A, t_0$ and $0 < k < 1$, and
  obtain
  \begin{displaymath}
    y(t) = \sn( \sqrt{I}(t-t_0);k), \qquad \text{where} \quad k = \frac{\omega}{\sqrt{I}},
  \end{displaymath}
which yields 
\begin{displaymath}
  \theta(t) = 2 \arcsin \left( \sn( \sqrt{I}(t-t_0);k) \right),
\end{displaymath}
for any $t_0 \in \R$. Given a particular solution $\theta(t)$, we can
compute $t_0$ from the equation $\theta(t_0) = 0$. This case
corresponds to circulating orbits, where the pendulum has enough
energy to go over the top. Note that we need to keep changing the
branches of $\arcsin$ in order to keep $\theta$ increasing or
decreasing.

\end{description}
  
\begin{lem}     \label{lem:int}
  Let $t \mapsto \theta(t)$ be a solution to $\ddot{\theta} + \omega^2
  \sin \theta = 0$, for $0 \leq t \leq \ell$. Suppose that
  \begin{equation}       \label{eq:int}
    \int_0^\ell \sin \frac{\theta(t)}{2} \: dt = \int_0^\ell \cos
    \frac{\theta(t)}{2} \: dt = 0.
  \end{equation}
  Then:

  \begin{enumerate}

  \item[(a)] $I > \omega^2$.

  \item[(b)] $\ell \sqrt{I}$ is an integer multiple of the period of
    $\: \sn( \: \cdot \:; k)$, with $k = \omega/\sqrt{I}$.

  \item[(c)] Furthermore,
    \begin{displaymath}
      \int_0^\ell \sin \theta(t) \: dt = 0.
    \end{displaymath}
  \end{enumerate}
\end{lem}

\begin{proof}
  The assumption \eqref{eq:int} implies that only Cases 2 and 4 above
  are possible. If
    \begin{displaymath}
      \sin \frac{\theta(t)}{2} = k \: \sn (\omega(t-t_0);k)
    \end{displaymath}
    as in Case 2, then by Proposition~\ref{prop:elliptic}(b),
    $\abs{\sin \theta(t)/2} \leq k$, so $\abs{\cos \theta(t)/2} \geq
    \sqrt{1-k^2}$, contradicting the second part of
    \eqref{eq:int}. Therefore, only Case 4 is possible, which implies
    $I > \omega^2$, proving (a). Furthermore, we know that
    \begin{displaymath}
      \sin \frac{\theta(t)}{2} = \sn \left(\sqrt{I}(t-t_0);k \right),
    \end{displaymath}
    where $k = \omega/\sqrt{I}$. Thus
    \begin{align*}
      \cos \frac{\theta(t)}{2} & = \pm \sqrt{1- \sin^2
        \frac{\theta(t)}{2}} \\
      & = \pm \sqrt{1 - \sn^2 \left(\sqrt{I}(t-t_0);k \right)} \\
      & = \pm \cn \left(\sqrt{I}(t-t_0);k \right).
    \end{align*}
    Since
    \begin{displaymath}
      \int_0^\ell \sin \frac{\theta(t)}{2} \: dt = \frac{1}{\sqrt{I}} 
      \int_{-\sqrt{I} \: t_0}^{\sqrt{I}(\ell-t_0)} \sn(t;k) \: dt
    \end{displaymath}
    and
    \begin{displaymath}
      \int_0^\ell \cos \frac{\theta(t)}{2} \: dt = \pm \frac{1}{\sqrt{I}} 
      \int_{-\sqrt{I} \: t_0}^{\sqrt{I}(\ell-t_0)} \cn(t;k) \: dt,
    \end{displaymath}
    Corollary~\ref{cor:elliptic} and \eqref{eq:int} imply that $\ell
    \sqrt{I}$ is an integer multiple of $4K$, where $4K$ is the period
    of $\sn(\: \cdot \:;k)$ and $\cn(\: \cdot \:;k)$. This proves (b).

    Finally, again using Corollary~\ref{prop:elliptic}, we obtain (c):
    \begin{align*}
      \int_0^\ell \sin \theta(t) \: dt & = \pm 2\int_0^\ell
      \sn(\sqrt{I}(t-t_0);k) \: \cn (\sqrt{I}(t-t_0);k) \: dt \\
      & = \pm \frac{2}{\sqrt{I}} \int_{-\sqrt{I} \:
        t_0}^{\sqrt{I}(\ell-t_0)} \sn (t;k) \: \cn(t;k) \: dt \\
      & = 0. \qedhere
    \end{align*}
\end{proof}

\section{Proof of Theorem A}
\label{sec:proofA}

Let $\Phi$ be a special contact Anosov flow on a closed 3-manifold
$M$. Then
\begin{equation}  \label{eq:brackets}
  [X,Y] = Y, \quad [Y,Z] = X, \quad \text{and} \quad [Z,X] = Z,
\end{equation}
where $Y \in E^{ss}$, $Z \in E^{uu}$ are $C^\infty$ vector fields and
$(X,Y,Z)$ is an orthonormal frame with respect to a fixed Riemannian
metric on $M$. On each fiber of $T^\ast M$ we introduce the
coordinates $(P_X,P_Y,P_Z)$, where $P_X, P_Y, P_Z$ are the momentum
functions of $X, Y, Z$, respectively. The subriemannian Hamiltonian is
\begin{displaymath}
  H = \frac{1}{2} (P_Y^2 + P_Z^2).
\end{displaymath}
The subriemannian geodesic equations in these coordinates are
\begin{align*}
  \dot{x} & = P_Y Y_x + P_Z Z_x \\
  \dot{P}_X & = \{ P_X, H \} \\
  \dot{P}_Y & = \{ P_Y, H \} \\
  \dot{P}_Z & = \{ P_Z, H \}.
\end{align*}
Using \eqref{eq:brackets} it is not hard to see that
\begin{displaymath}
  \{ P_X, H\} = P_Z^2 - P_Y^2, \quad \{ P_Y, H \} = P_X P_Z, \quad
  \text{and} \quad \{ P_Z,H \} = - P_X P_Y.
\end{displaymath}
Thus the geodesic equations are
\begin{align*}
  \dot{x} & = P_Y Y_x + P_Z Z_x \\
  \dot{P}_X & = P_Z^2 - P_Y^2 \\
  \dot{P}_Y & = P_X P_Z \\
  \dot{P}_Z & = -P_X P_Y.
\end{align*}
Since $H$ is an integral of motion, we will consider only solutions
lying on the level set $\ms{L}$ defined by $P_Y^2 + P_Z^2 = 1$. Each
fiber of $\ms{L}$ is topologically a cylinder $S^1 \times \R$. We set
\begin{displaymath}
  P_Y = \cos \theta, \qquad P_Z = \sin \theta,
\end{displaymath}
so that on $\ms{L}$ we have local coordinates $(x;\theta,P_X)$, where
$x \in M$. In these coordinates the subriemannian Hamiltonian
equations become
\begin{align*}
  \dot{x} & = (\cos \theta) Y + (\sin \theta) Z \\
  \dot{\theta} & = -P_X \\
  \dot{P}_X & = \sin^2 \theta - \cos^2 \theta.
\end{align*}
Therefore, $\theta$ satisfies 
\begin{displaymath}
  \ddot{\theta} - \cos 2\theta = 0.
\end{displaymath}
Let $x \in M$ be arbitrary and pick a flowbox $U$ for $X$ containing
$x$. Let $\tau > 0$ be small enough so that $y = f_\tau(x)$ lies in
$U$ and let $\gamma$ be a unit speed subriemannian geodesic such that
$\gamma(0) = x$ and $\gamma(\ell) = y$, where $\ell =
d_{su}(x,y)$. Then
\begin{displaymath}
  \dot{\gamma}(t) = \cos \theta(t) Y_{\gamma(t)} + \sin \theta(t) Z_{\gamma(t)},
\end{displaymath}
where $\theta$ satisfies the harmonic oscillator equation
$\ddot{\theta} - \cos 2\theta = 0$. Observe that $\gamma$ is
$su$-balanced if
\begin{displaymath}
  \int_0^\ell \cos^2 \theta(t) \: dt = \int_0^\ell \sin^2 \theta(t) \: dt,
\end{displaymath}
which is equivalent to 
\begin{displaymath}
  \int_0^\ell \cos 2 \theta(t) \: dt = 0.
\end{displaymath}
The substitution $\phi = 2\theta - \dfrac{\pi}{2}$ in the
$\theta$-equation converts it to the standard form
\begin{displaymath}
  \ddot{\phi} + 2 \sin \phi = 0,
\end{displaymath}
with $\omega = \sqrt{2}$.

Now consider the orbit space $\Sigma = U/\! \sim$ of the Anosov flow
in $U$, where $z \sim w$ if $w$ is on the orbit of $z$. It is clear
that $\Sigma$ is diffeomorphic to an open set in $\R^2$. Denote by
$\tilde{\gamma}$ by the natural projection of $\gamma$ to $\Sigma$ and
by $\tilde{Y}$, and $\tilde{Z}$ the projections of $Y$ and $Z$ along
$\gamma$ to $\Sigma$. Then relative to the frame
$(\tilde{Y},\tilde{Z})$ along $\tilde{\gamma}$, we have
\begin{displaymath}
  \dot{\tilde{\gamma}}(t) = (\cos \theta(t), \sin \theta(t)),
\end{displaymath}
for $0 \leq t \leq \ell$. Since $\tilde{\gamma}(\ell) =
\tilde{\gamma}(0)$, it follows that 
\begin{displaymath}
  \int_0^\ell \cos \theta(t) \: dt = \int_0^\ell \sin \theta(t) \: dt
  = 0.
\end{displaymath}
Using $\theta(t) = \dfrac{1}{2} \phi(t) + \dfrac{\pi}{4}$ and
substituting, we obtain 
\begin{equation}  \label{eq:phi-int}
  \int_0^\ell \cos \frac{\phi(t)}{2} \: dt = \int_0^\ell \sin \frac{\phi(t)}{2} \: dt = 0.
\end{equation}
Lemma~\ref{lem:int} yields
\begin{displaymath}
  \int_0^\ell \sin \phi(t) \: dt = 0.
\end{displaymath}
But $\sin \phi = -\cos 2 \theta = \sin^2 \theta - \cos^2 \theta$, so 
\begin{displaymath}
  \int_0^\ell \sin^2 \theta(t) \: dt = \int_0^\ell \cos^2 \theta(t) \: dt
  = 0,
\end{displaymath}
proving that $\gamma$ is an $su$-balanced geodesic.

\section{Proof of Theorem B}
\label{sec:proofB}

Assume $\Phi = \{ f_t \}$ is the geodesic flow of a closed Riemannian
manifold $N$ with constant negative curvature; $\Phi$ is defined on
$M$, the unit tangent bundle of $N$. Without loss we can assume that
\begin{displaymath}
  \norm{Tf_t(v)} = e^{-t} \norm{v} \qquad \text{and} \qquad \norm{Tf_t(w)} = e^t \norm{w},
\end{displaymath}
for all $v \in E^{ss}$, $w \in E^{uu}$ and all real $t$.

Let $p \in M$ and $\tau > 0$, and set $q = f_\tau(p)$. Let $\gamma :
[0,\ell] \to M$ be a unit speed $su$-subriemannian geodesic from $p$
to $q$ so that $\abs{\gamma} = \ell = d_{su}(p,q)$. We assume $\tau$
is sufficiently small so that the closed curve consisting of $\gamma$
and the arc of $X$-orbit $\{ f_t(p) : 0 \leq t \leq \tau \}$ is
homotopically trivial. Clearly, there exists a uniform $\delta > 0$
such that this is true for every two points on the same $\Phi$-orbit
which are $\delta$-close (with respect to the Riemannian distance) and
for every $su$-subriemannian geodesic connecting them.

We claim that
\begin{equation}    \label{eq:t}
  d_{su}(f_t(p),f_t(q)) = d_{su}(p,q),
\end{equation}
for all $t \in \R$.

To prove this, first denote by $\tilde{p}, \tilde{q}$ the lifts of $p,
q$ lying in the same fundamental domain of the universal Riemannian
covering $\pi : \tilde{M} \to M$. Let $t \in \R$ be arbitrary and let
$x = f_t(p)$ and $y = f_t(q)$. Denote by $\tilde{x}, \tilde{y}$ the
corresponding lifts to the universal covering, so that $\tilde{x} =
\tilde{f}_t(\tilde{p})$ and $\tilde{y} = \tilde{f}_t(\tilde{q})$.

By Lemma~\ref{lem:hyperbolic} there exists an isometry $F$ of
$\tilde{M}$ such that $F(\tilde{p}) = \tilde{x}$ and $TF(\tilde{X}) =
\tilde{X}$, where $\tilde{X}$ denotes the lift of $X$ to
$\tilde{M}$. Denote by $\tilde{\gamma}$ the lift of $\gamma$. Then $F
\circ \tilde{\gamma}$ connects $\tilde{x}$ and $\tilde{y}$ and has the
same length as $\tilde{\gamma}$ as well as $\gamma$.

\begin{figure}[htbp]
\centerline{
	\psfrag{x}[][]{$x$}
	\psfrag{y}[][]{$y$}
	\psfrag{p}[][]{$p$}
	\psfrag{q}[][]{$q$}
        \psfrag{g}[][]{$\gamma$}
        \psfrag{pt}[][]{$\tilde{p}$}
        \psfrag{qt}[][]{$\tilde{q}$}
        \psfrag{xt}[][]{$\tilde{x}$}
        \psfrag{yt}[][]{$\tilde{y}$}
        \psfrag{gt}[][]{$\tilde{\gamma}$}
        \psfrag{gy}[][]{$\hat{\gamma} = \pi \circ F \circ
          \tilde{\gamma}$}
        \psfrag{Fgt}[][]{$F \circ \tilde{\gamma}$}
        \psfrag{pi}[][]{$\pi$}
\includegraphics[width=0.4\hsize]{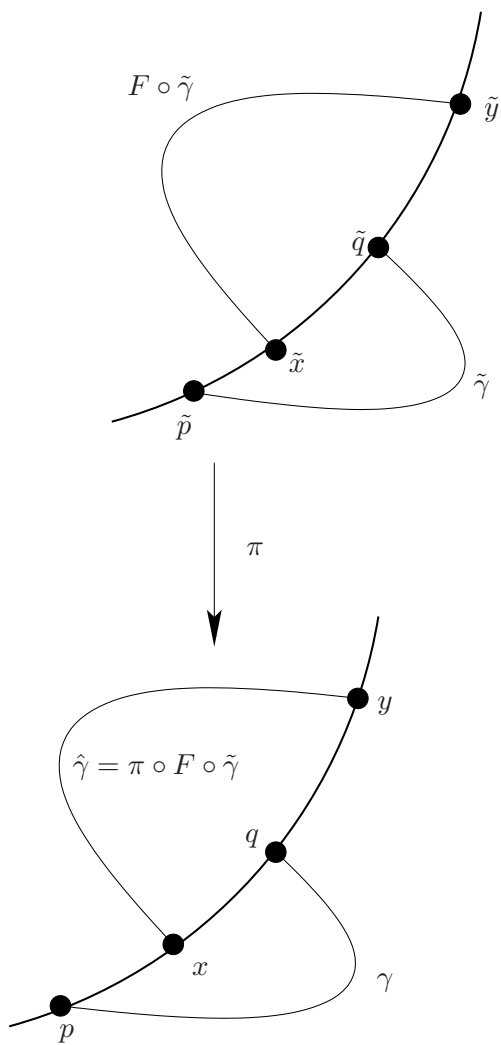}}
\caption{Proof of Theorem B.}
\label{fig:proofB}
\end{figure}

We claim that its projection $\hat{\gamma} = \pi \circ F \circ
\tilde{\gamma}$ is an $su$-geodesic connecting $x$ and $y$. See
Figure~\ref{fig:proofB}. First observe that since $F$ preserves
$\tilde{X}$, and $\tilde{E}^{su}$ and $\tilde{X}$ are orthogonal,
$\hat{\gamma}$ is an $su$-path. Suppose that it is not minimal. Then
there exists a horizontal path $c$ joining $x$ and $y$ such that
$\abs{c} < \abs{\hat{\gamma}}$. This path lifts to a path $\tilde{c}$
in $\tilde{M}$ connecting $\tilde{x}$ and $\tilde{y}$, and tangent to
the lift of $E^{su}$; $\tilde{c}$ has the same length as $c$. It
follows that $\pi \circ F^{-1} \circ \tilde{c}$ is a horizontal path
joining $p$ and $q$ whose length is
\begin{displaymath}
  \abs{\pi \circ F^{-1} \circ \tilde{c}} = \abs{c} <
  \abs{\hat{\gamma}} = \abs{\gamma},
\end{displaymath}
contradicting the assumption that $\gamma$ is a geodesic. This
proves \eqref{eq:t}.

Now let us show that $\gamma$ is $su$-balanced. Write $\dot{\gamma}(t)
= w_s(t) + w_u(t)$, with $w_s(t) \in E^{ss}$ and $w_u(t) \in E^{uu}$,
for all $0 \leq t \leq \ell$. Then:
\begin{align}
  \left. \frac{d}{dr} \right|_{r=0} \abs{f_r \circ \gamma} & = 
     \left. \frac{d}{dr} \right|_{r=0} \int_0^\ell 
     \norm{Tf_r(\dot{\gamma}(t))} \: dt \nonumber \\
     & =  \left. \frac{d}{dr} \right|_{r=0} \int_0^\ell 
     \norm{Tf_r(w_s(t)) + Tf_r(w_u(t))} \: dt \nonumber \\
     & =  \left. \frac{d}{dr} \right|_{r=0} \int_0^\ell
     \sqrt{\norm{Tf_r(w_s(t))}^2 + \norm{Tf_r(w_u(t))}^2} \: dt \nonumber \\
     & = \left. \frac{d}{dr} \right|_{r=0} \int_0^\ell
     \sqrt{e^{-2r} \norm{w_s(t)}^2 + e^{2r} \norm{w_u(t)}^2} \: dt
     \nonumber \\
     & = 2 \int_0^\ell \left\{ \norm{w_u(t)}^2 - 
         \norm{w_s(t)}^2 \right\} \: dt \nonumber \\
     & = 2 \ell \{\mathcal{E}_u(\gamma) - \mathcal{E}_s(\gamma) \}.   \label{eq:diff}
\end{align}
Here we used the fact that $\gamma$ is unit-speed and $w_s(t)$ and
$w_u(t)$ are orthogonal. If $\gamma$ is not $su$-balanced, i.e.,
$\mathcal{E}_s(\gamma) \neq \mathcal{E}_u(\gamma)$, then by \eqref{eq:diff}, there exists
$t \neq 0$ such that $\abs{f_t \circ \gamma} < \abs{\gamma}$, which
yields
\begin{displaymath}
  d_{su}(f_t(p),f_t(q)) \leq \abs{f_t \circ \gamma} < \abs{\gamma} = d_{su}(p,q),
\end{displaymath}
contradicting \eqref{eq:t}. This completes the proof.

\section{Open questions}
\label{sec:questions}

Virtually any question involving subriemannian geometry defined by the stable and unstable bundles of a partially hyperbolic dynamical system or an Anosov flow is open. Here we list only a few.

\begin{enumerate}

\item If $\Phi$ is an \emph{arbitrary} contact Anosov flow on a 3-manifold, what can be said about $\mathcal{E}_s(\gamma)$ and $\mathcal{E}_u(\gamma)$, for an $su$-subriemannian geodesic $\gamma$ whose endpoints lie on an orbit of the flow? Are there lower and upper bounds for $\mathcal{E}_s(\gamma)/\mathcal{E}_u(\gamma)$ independent of $\gamma$? The same question can be asked about arbitrary contact Anosov flows in any odd dimension. 

\item What if $\Phi$ is the geodesic flow of a closed Riemannian manifold with \emph{variable} negative sectional curvature?

\item Suppose $\Phi$ is a transitve Anosov flow such that $E^{su}$ has the accessibility property but is not necessarily contact. Do $su$-subriemannian geodesics exist? If so, how close are they to being $su$-balanced (when their endpoints lie on the same orbit of $\Phi$)?

\item The same question can be asked for an accessible partially hyperbolic diffeomorphism.

\end{enumerate}


\bibliographystyle{amsalpha} 

\begin{thebibliography}{HPS77}

\bibitem[Ano67]{anosov+67}
Dimitri~V. Anosov, \emph{Geodesic flows on closed {R}iemannian manifolds of
  negative curvature}, Proc. {S}teklov {M}ath. {I}nst. \textbf{90} (1967),
  {AMS} {T}ranslations (1969).

\bibitem[AS67]{anosov+sinai+67}
Dimitri~V. Anosov and Yakov~G. Sinai, \emph{Some smooth ergodic systems}, Russ.
  {M}ath. {S}urveys \textbf{22} (1967), 103--167.

\bibitem[Boo03]{boothby+2003}
William~M. Boothby, \emph{An introduction to differentiable manifolds and
  {R}iemannian geometry}, Academic Press, 2003.

\bibitem[Ebe73]{eberlein+73}
P.~Eberlein, \emph{When is a geodesic flow of {A}nosov type? {I}, {II}},
  Journal of Differential Geom. \textbf{8} (1973), pp. 437--463, ibid. 8
  (1973), pp. 565--577.

\bibitem[FH13]{foulon+hass+13}
Patrick Foulon and Boris Hasselblatt, \emph{Contact {A}nosov flows on
  hyperbolic 3--manifolds}, Geom. Topol. \textbf{17} (2013), 1225--1252.

\bibitem[Ghy87]{ghys+87}
Etienne Ghys, \emph{Flots d'{A}nosov dont les feuilletages stables sont
  diff{\'e}rentiables}, Ann. Sci. \'Ecole Norm. Sup. (4) \textbf{20} (1987),
  no.~2, 251--270.

\bibitem[HK90]{hurder+katok+90}
Steven Hurder and Anatole Katok, \emph{Differentiability, rigidity and
  {G}odbillon-{V}ey classes for {A}nosov flows}, Publ. Math. IHES \textbf{72}
  (1990), no.~1, 5--61.

\bibitem[HP75]{hp+75}
Morris~W. Hirsch and Charles~C. Pugh, \emph{Smoothness of horocycle
  foliations}, Journal of Differential Geom. \textbf{10} (1975), 225--238.

\bibitem[HPS77]{hps77}
Morris~W. Hirsch, Charles~C. Pugh, and Michael Shub, \emph{Invariant
  manifolds}, Lecture {N}otes in {M}athematics, vol. 583, Springer-Verlag,
  Berlin-New York, 1977.

\bibitem[Liv04]{liverani+contact+04}
Carlangelo Liverani, \emph{On contact {A}nosov flows}, Ann. of Math. (2)
  \textbf{159} (2004), no.~3, 1275--1312.

\bibitem[Mey01]{meyer+amm+01}
Kenneth~R. Meyer, \emph{Jacobi elliptic functions from a dynamical systems
  point of view}, Amer. Math. Montly \textbf{108} (2001), 729--737.

\bibitem[Mon02]{mont02}
Richard Montgomery, \emph{A tour of subriemannian geometries, their geodesics
  and applications}, Mathematical {S}urveys and {M}onographs, vol.~91, AMS,
  2002.

\bibitem[MS99]{duff+solomon+sympl+99}
Dusa McDuff and Dietmar Salamon, \emph{Introduction to symplectic topology},
  second ed., Oxford Mathematical Monographs, Oxford University Press, 1999.

\bibitem[Par86]{parry86}
William Parry, \emph{Synchronisation of canonical measures for hyperbolic
  attractors}, Comm. Math. Phys. \textbf{106} (1986), no.~2, 267--275.

\bibitem[Pat99]{paternain+geodesic+99}
Gabriel~P. Paternain, \emph{Geodesic flows}, Progress in {M}athematics, vol.
  190, Birkh\"auser, 1999.

\bibitem[Pla72]{plante72}
Joseph Plante, \emph{Anosov flows}, Amer. J. of Math. \textbf{94} (1972),
  729--754.

\bibitem[PS04]{ps+04}
Charles~C. Pugh and Michael Shub, \emph{Stable ergodicity}, Bull. Amer. Math.
  Soc. \textbf{41} (2004), 1--41.

\bibitem[Sim97]{sns+97}
Slobodan~N. Simi\'c, \emph{Codimension one {A}nosov flows and a conjecture of
  {V}erjovsky}, Ergod. Th. Dynam. Syst. \textbf{17} (1997), 1211--1231.

\bibitem[Sim10]{sns+pams+10}
Slobodan~N. Simi\'{c}, \emph{A lower bound on the subriemannian distance for
  {H}{\"o}lder distributions}, Proc. Amer. Math. Soc. \textbf{138} (2010),
  no.~9, 3293--3299.

\end{thebibliography}

\end{document}